\documentclass{amsart}
\usepackage{amssymb}
\newcommand{\PM}{{\mathbf {PM}}}
\newcommand{\UM}{{\mathbf  {UM}}}
\newcommand{\UN}{{\mathbf  {UN}}}
\newcommand{\EE}{\mathbb{EE}}
\newcommand{\PP}{\mathbb{P}}
\newcommand{\X}{\mathbf{X}}

\newcommand{\spli}{{\operatorname{\mathsf   {split}}}}

\newcommand{\suc}{{\operatorname{\mathsf    {succ}}}}

\newcommand{\forces}{\Vdash}
\newcommand{\ZFCa}{{\operatorname{\mathsf {ZFC}}}}
\newcommand{\CH}{\operatorname{\mathsf {CH}}}
\newcommand{\reals}{{\mathbb R}}

\newcommand{\rest}{{\mathord{\restriction}}}

\newcommand{\unif}{\operatorname{\mathsf  {non}}}

\newcommand{\dom}{{\operatorname{\mathsf {dom}}}}
\newcommand{\supp}{{\operatorname{\mathsf {supp}}}}
\newcommand{\cl}{{\operatorname{\mathsf {cl}}}}
\newcommand{\ot}{{\operatorname{\mathsf {ot}}}}

\newcommand{\N}{{\mathcal N}}
\newcommand{\M}{{\mathcal M}}
\newcommand{\V}{{\mathbf V}}
\newcommand{\<}{\langle}
\renewcommand{\>}{\rangle}

\newcommand{\thinks}{\models}

\newcommand{\lft}[2]{\mathopen\ifcase#1{}\oo\or
                        \big#2\or\Big#2\else\oo\fi} 
\newcommand{\rgt}[2]{\mathclose\ifcase#1{}\oo\or
                        \big#2\or\Big#2\else\oo\fi}

\theoremstyle{plain}
\newtheorem{theorem}{Theorem}
\theoremstyle{plain}
\newtheorem{lemma}[theorem]{Lemma}

\newtheorem{definition}[theorem]{Definition}

\begin{document}
\title{Perfectly meager sets and universally null sets}
\author{Tomek Bartoszynski}
\address{Department of Mathematics and Computer Science\\
Boise State University\\
Boise, Idaho 83725 U.S.A.}
\thanks{The first author was  partially supported by 
NSF grant DMS 9971282} 
\email{tomek@math.idbsu.edu, http://math.idbsu.edu/\char 126 tomek}
\author{Saharon Shelah}
\thanks{The second author was partially supported by Israel Science
   Foundation. Publication 732}
\address{Department of Mathematics\\
Hebrew University\\
Jerusalem, Israel}
\email{shelah@math.huji.ac.il, http://math.rutgers.edu/\char 126 shelah/}
\begin{abstract}
We will show that
there is no $\ZFCa$ example of a set distinguishing between
universally null and perfectly meager sets.
\end{abstract}
\maketitle

\section{Introduction}
Consider the following three families of sets of reals:
\begin{definition}
  Let $X \subseteq \reals$.
  \begin{enumerate}
  \item $X$ is perfectly meager if for every perfect set $P \subseteq
    \reals$, $P \cap X$ is meager in $P$.
  \item $X$ is universally meager if every Borel isomorphic image of
    $X$ is meager,
  \item $X$ is universal null if every Borel isomorphic image of
    $X$ has Lebesgue measure zero.
  \end{enumerate}
Let $\PM$, $\UM$ and $\UN$ denote these families respectively.
\end{definition}
One gets an equivalent definition of $\UN$ by replacing
``Borel isomorphic'' by ``homeomorphic'', but this is not the case
with $\UM$.

Let $\M$ and $\N$ denote the $\sigma $-ideals of meager and of measure
zero subsets of the reals, respectively. 

For a $\sigma $-ideal $  {\mathcal J} \subseteq P(\reals)$ let 
$$\unif({\mathcal J})=\min\{|X|: X \subseteq \reals \ \&\ X \not
\subseteq {\mathcal J} \}.$$

There are many $\ZFCa$ examples of uncountable sets that are in $\UM \cap
\UN$. These include $\omega_1 \omega_1^\star$-gaps, a selector from the
constituents of a non-Borel ${\boldsymbol \Pi}^1_1$ set, etc. (see \cite{Mil84Spe})
All these sets have size $\boldsymbol\aleph_1 $, since Miller
\cite{Mil83Map} showed that, consistently,
no set of size $2^{\boldsymbol\aleph_0} $ is in $\UM \cup \UN$.

Grzegorek found  other constructions in $\ZFCa$ that produce sets of
(consistently) different sizes.
\begin{theorem}[Grzegorek, \cite{Grz80Sol}]
  \begin{enumerate}
  \item There exists a set $X \in \UN$ such that $|X|=\unif(\N)$,
  \item There exists a set $X \in \UM$ such that $|X|=\unif(\M)$.
  \end{enumerate}
\end{theorem}

The problem
whether the equality $\UM=\UN$ is consistent is open. However, both inclusions are
consistent with $\ZFCa$; $\UM \subseteq \UN$ holds in a model obtained by
adding ${\boldsymbol\aleph}_2 $ Cohen reals, and $\UN \subseteq \UM$ holds
in a model obtained by adding $ {\boldsymbol\aleph}_2 $ random reals
(side-by-side) (see \cite{Mil84Spe}, \cite{Mil83Map}).

In this paper we investigate the connection between families $\UN$ and
$\PM$, and show that both inclusions $\PM \subseteq \UN$ and $\UN
\subseteq \PM$ are consistent with $\ZFCa$ as well.
Observe that trivially $\UM \subseteq \PM$, thus we only need to check
that $\PM \subseteq
\UN$ is consistent. 
Recall that $\PM \neq \UM$ is consistent (\cite{Sie}) as well as
$\PM=\UM$ (\cite{BaPerf00}).

We will show that:
\begin{theorem}
  It is consistent with $\ZFCa$ that
$$\PM \subseteq [\reals]^{\leq \boldsymbol\aleph_1 } \subseteq \UN.$$
\end{theorem}

\section{Forcing}
  Suppose that $X \subseteq 2^\omega $ is a perfectly meager set in
  $\V$. Let $\widetilde{P}$ be a fixed  closed subset of $2^\omega \times
  2^\omega $ which is universal for perfect sets in $2^\omega $. In
  other words, for ever perfect set $P \subseteq 2^\omega $ there
  exists $x$ such that  $P=(\widetilde{P})_x$.
Since $X$ is perfectly meager, we can find a sets $\widetilde{Q^n}
\subseteq 2^\omega \times 2^\omega $ such that  for every $x \in
2^\omega$, and $n \in \omega $,
\begin{enumerate}
\item $(\widetilde{Q^n})_x$ is a closed nowhere dense subset of $(\widetilde{P})_x$,
\item $X \cap (\widetilde{P})_x  \subseteq (\bigcup_{n \in
    \omega}\widetilde{Q^n})_x$. 
\end{enumerate}
Clearly, the set $\bigcup_{n \in \omega} \widetilde{Q^n}$ witnesses
that $X \in \PM$ since 
$$X \subseteq 2^\omega \setminus \bigcup_{x \in 2^\omega } (\widetilde{P}
\setminus \bigcup_{n \in \omega} \widetilde{Q^n})_x .$$

Note that the last inclusion makes sense even if $X$ is not a subset
of 
$\V$.
Suppose that $\V' \subseteq \V$ and $X \subseteq \V$ is a set of
reals. We will say $\V' \thinks X \in \PM$ if
there exists a family $\{\widetilde{Q^n}: n \in \omega\} \in \V'$ such
that  $ X \cap (\widetilde{P})_x  \subseteq (\bigcup_{n \in
    \omega}\widetilde{Q^n})_x$ for every real $x \in \V'$.

The property of being perfectly meager is not absolute so whether $X$
is perfectly meager in $\V'$ has no bearing onto whether $X$ is
perfectly meager in $\V$. For example, if $x \in \V$ is a Cohen real
over $\V'$ then the set $\{x\}$ is perfectly meager in $\V$ but not in
$\V'$. 
\begin{lemma}
  Let $\<{\mathcal P}_\alpha, \dot{{\mathcal Q}}_\alpha :
  \alpha<\omega_2\}$ be a countable support iteration  of proper
  forcing notions over $\V \thinks \CH$. Suppose that $X \subseteq
  \V^{{\mathcal P}_{\omega_2}} \cap \reals$ is a perfectly meager
  set. Then there exists an
  $\omega_1$-club $C \subseteq \omega_2$ such that  for every $\alpha
  \in C$,
$$\V^{{\mathcal P}_\alpha} \thinks X \in \PM.$$
\end{lemma}
\begin{proof}
  Let $\{\widetilde{Q^n}: n \in \omega\} \in \V^{{\mathcal
      P}_{\omega_2}}$ be a family witnessing that $X$ is perfectly
  meager. Let $C$ consist of those ordinals of cofinality $\omega_1$
  that for every $n$, $\widetilde{Q^n} \cap \lft1((2^\omega \cap
  \V^{{\mathcal P}_\alpha})\times 2^\omega \rgt1) \in \V^{{\mathcal
      P}_\alpha}$. The usual argument involving Skolem-L\"owenheim
  theorem shows that $C$ has the required property.
\end{proof}

Our objective is to find a set of general conditions on a forcing
notion $\PP$ such that the countable support iteration of $\PP$ of
length $\omega_2$ produces a model where $\PM \subseteq [\reals]^{\leq
  \boldsymbol\aleph_1} \subseteq \UN$.
These conditions are sufficient for the class of forcing notions
defined using norms \cite{RoSh470}.

These conditions are the following:

\begin{enumerate}
\item $\V^{\PP} \thinks \V \cap 2^\omega \in \N$,
\item $\V^{\PP} \thinks \V \cap 2^\omega \not\in \M$,
\item $\PP$ is $\omega^\omega $-bounding, that is $\omega^\omega \cap
  \V$ is a dominating family in $\omega^\omega \cap \V^{\PP}$,
\item $\PP$ adds a real $x_\PP \in 2^\omega $ such that $\V \thinks
  \{x_\PP\} \not\in \PM$.
\item $\PP$ generic real is minimal, that is, if $g$ is $\PP$-generic
  over $\V$ and $x \in \V[g] \cap 2^\omega $ then $x \in \V$ or $g
  \in \V[x]$.
\end{enumerate}

Condition (1) is necessary to make all sets of size
${\boldsymbol\aleph_1}$ universally null, and condition (2) is
necessary to avoid making  all ${\boldsymbol\aleph_1} $ sets perfectly
meager. 
Recall that (2) and (3) together are essentially equivalent to 
$$\V^{\PP} \thinks \V \cap \M \text{ is cofinal in } \M.$$

For the forcing notions $\PP$ that we have in mind the following property
holds: for every real $x \in \V^{\PP}$ there exists a continuous
function $f \in \V$ such that  $x=f(x_G)$, where $x_G$ is a generic real.

Condition (5), guarantees that in the above context $f$ can
be chosen to  be a homeomorphism. In particular, if $X$ is a
set of reals of size $\boldsymbol\aleph_2$ then $X$ will contain a
homeomorphic image of a sequence of generic reals.

The following forcing notion appeared in \cite{GJS92}, it is similar (but not
identical) to the infinitely equal real forcing from \cite{Mil81Som}.

For a tree $p$ and $t \in p$,
let $\suc_p(t)$ be the set of
all immediate successors of $t$ in $p$, $p_t=\{v \in p: t
\subseteq v \text{ or } v \subseteq t\}$ the subtree of $p$ determined
by $t$, $p\rest n$ the $n$-th level of $p$, and let $[p]$ be the set of
branches of $p$. By identifying $s \in \omega^{<\omega}$ with the
full-branching tree having root $s$, we can also denote
$[s]=\{f \in \omega^\omega:s \subseteq f\}$.

Fix a strictly increasing function $f \in \omega^\omega $ and let
$\X=\prod_{n \in \omega} f(n)$. Note that $\X$ is a Polish space
homeomorphic to $2^\omega $. For technical reasons we require that 
$f(n)=2^{\tilde{f}(n)}$ for $n \in \omega $.

Let $\EE$ be the following forcing notion:
\begin{multline*}
p \in \EE \iff p \subseteq \bigcup_{n \in \omega}\bigcup_{j<n} \prod f(j) \text{ is a
  perfect tree } \& \\ 
\forall s \in p \ \exists t \in p \
\lft2(s \subseteq t \ \& \ \suc_p(t)=f(|p|)\rgt2).
\end{multline*} 
For $p, q \in \EE$, $p \geq q$ if $p \subseteq q$.
Without loss of generality we can assume that  $|\suc_p(s)|=1$ or
$\suc_p(s)=f(|p|)$ for all $p\in \EE$ and $s \in
p$.
Conditions of this type form a dense subset of $\EE$.

Let      
$$ \spli(p) = \{s\in p : |\suc_p(s)|>1\} = \bigcup_{n \in \omega} \spli_n(p),$$
where
$\spli_n(p)=\left\{s \in \spli(p): \lft2|\lft1\{t \subsetneq s: t \in
\spli(p)\rgt1\}\rgt2| =n\right\}$.

For $p,q\in \EE$, $n \in \omega$,   we let 
$$p \geq_n q \iff p \geq q \ \& \ \spli_n(q)=\spli_n(p).$$

\begin{lemma}[\cite{GJS92}]\label{old}
  \begin{enumerate}
  \item $\EE$ satisfies Axiom A, so it is proper,
  \item $\V^{\EE} \thinks \V \cap 2^\omega \in \N$,
\item $\V^{\EE} \thinks \V \cap 2^\omega \not\in \M$,
    \item for every maximal antichain ${\mathcal A} \subseteq \EE$,
    $p\in \EE$, and
    $n \in \omega $ there exists $q \geq_n p$ such that $\{r \in
    {\mathcal A} : r \text{ is compatible with } q\}$ is finite.
  \item for every family of maximal antichains $\{{\mathcal A}_n: n
    \in \omega \}$ and
    $p\in \EE$
    there exists $q \geq p$ such that for every $n$, 
$\{r \in
    {\mathcal A}_n : r \text{ is compatible with } q\}$ is finite.
  \item $\EE$ is $\omega^\omega $ bounding,
\item $\V^{\EE} \thinks \V \cap \M \text{ is cofinal in } \M$. $\qed$
  \end{enumerate}
\end{lemma}

Note that for $p \in \EE$ 
the set $[p]$ is a compact subset of $\X=\prod_n f(n)$. 
Moreover, there is a canonical isomorphism between $[p]$ and $2^\omega
$ defined as follows:

For every $n$ let $\{s^n_0, \dots, s^n_{f(n)}\}$ be a fixed
enumeration of 0-1 sequences of length $\tilde{f}(n)$ (recall that
$f(n)=2^{\tilde{f}(n)}$). 
Define $F: [p] \longrightarrow 2^\omega $ as
$$F(x)=s^{n_0}_{x(n_0+1)}{}^\frown s^{n_1}_{x(n_1+1)}{}^\frown \dots
,$$
where $n_0, n_1, \dots $ is the increasing enumeration of the set 
$\{n: x \rest n \in \spli(p)\}$.
\begin{lemma}\label{spec}
  Let $p \in \EE$ and suppose that $H \subseteq [p]$ is a
  meager set in $[p]$. For every $n \in \omega $ there
  exists $q \geq_n p$ such that $[q] \cap H=\emptyset$. In particular,
$$\forces_{\EE} \text{``}\V \thinks \{\dot{g}\} \not\in \PM\text{''} .$$
\end{lemma}
\begin{proof}
Let $H \subseteq [p]$ be a meager set, and let $n \in \omega $. Fix a descending sequence of
open sets $\<U_k: k \in \omega\>$ such that each $U_k$ is dense in
$[p]$ and $H \cap \bigcap_k U_k = \emptyset$.
By induction build a sequence $\<p_k: k \in \omega\>$ such that
$p_0=p$, and for every $k$,
\begin{enumerate}
\item $p_{k+1}\geq_{n+k+1} p_k \in \EE$,
\item $[p_{k+1}] \subseteq U_k$.
\end{enumerate}
Suppose that $p_k$ is given. For every $v \in \spli_{n+k+1}(p_{k})$
find $q_v \geq (p_{k})_v$ such that $[q_v] \subseteq U_k$. 
Let $p_{k+1}=\bigcup \{q_v: v \in \spli_{n+k+1}(p_k)\}$.
 Condition $q = \lim_k p_k$ has the required property. 

Suppose that $\{\widetilde{Q^n}: n \in \omega\} \in \V$ is a possible
witness that $\{\dot{g}\}$ is perfectly meager, and let $p \in
\EE$. Find $x \in \V$ such that  $[p]=(P)_x$ and let $q \geq p$ be such that
$[q] \cap \left(\bigcup_n \widetilde{Q^n}\right)_x=\emptyset$.
Clearly,
$$q \forces_{\EE} \{\dot{g}\} \in \bigcup_{x \in \V} \left(P \setminus
  \bigcup_n \widetilde{Q^n}\right)_x.$$
In particular,
$$q \forces_{\EE} \text{``$\V \models \{\dot{g}\} \not\in \PM$''}.$$
\end{proof}

\begin{lemma}\label{general}
  Suppose that $p \in \EE$ and $p \forces_{\EE} \dot{x} \in 2^\omega
  $. For every $n \in \omega $ there exists $q \geq_n p$ and a
  continuous function $F: [q] \longrightarrow 2^\omega $
  such that 
$$q \forces_{\EE} \dot{x}=F(\dot{g}),$$
where $\dot{g}$ is the canonical name for the generic real.

Moreover, we can require that for every $v \in \spli_n(q)$ and any
$x_1, x_2 \in [q_v]$, $F(x_1) \rest n = F(x_2) \rest n$.
\end{lemma}
\begin{proof}
The first part is a special case of a more general fact.
For $n \in \omega $ let ${\mathcal A}_n \subseteq \EE$ be a maximal
antichain below $p$ such that 
$$\forall r \in {\mathcal A}_n \ \exists s \in 2^n \ r \forces_{\EE}
\dot{x}\rest n =s.$$
Use \ref{old}(5) to find $q \geq p$ such that for every $n \in \omega
$,
$$\{r \in
    {\mathcal A}_n : r \text{ is compatible with } q\}$$ is finite.
Let ${\mathcal A}_n'=\{r \in
    {\mathcal A}_n : r \text{ is compatible with } q\}$.
Without loss of generality we can assume that $[q] \subseteq
\bigcup_{r \in {\mathcal A}_n'} [r]$.
It follows that $[r] \cap [q]$ is clopen in
$[q]$ for every $r \in {\mathcal A}_n'$.
Define $F:[q] \longrightarrow 2^\omega $ as $F(x)=y$ if for
every $n \in \omega $ there exists $r \in {\mathcal A}_n'$ such that
$x \in [r]$ and $r \forces_{\EE} \dot{x} \rest n = y\rest
n$.
It is easy to see that $F$ is a continuous function that has the
required properties. 

To show the second part we need to build $q$ in such a way that for
every $v \in \spli_n(q)$, there is $r \in {\mathcal A}_n'$ such that
$q_v \geq r$.
\end{proof}

\begin{lemma}\label{uff}
  Suppose that $p \in \EE$, $n \in \omega $ and $p \forces_{\EE} \dot{x} \in 2^\omega
  $. Let $F : [q] \longrightarrow 2^\omega $ be a continuous function
  such that $p \forces_{\EE} \dot{x}=F(\dot{g})$. 

There exists $q \geq p$ such that $ F \rest [q]$ is constant, 
or there exists $q \geq_n p$ such that $F \rest [q]$ is one-to-one. In
particular, the generic real is minimal.
\end{lemma}
\begin{proof}
Consider the following two cases.

{\sc Case 1} $p\not \forces_{\EE} \dot{x} \not\in
\V$. Let $x \in \V$ and $q \geq p$ be such that $q \forces_{\EE}
\dot{x}=x$. Clearly $F \rest [q]$ is constant with value $x$.

\bigskip 

{\sc Case 2} $p \forces_{\EE} \dot{x} \not \in \V$.

Build by induction a sequence of conditions $\<p_k:k\in \omega \>$
such that $p_0=p$ and
for every $k$,
\begin{enumerate}
\item $p_{k+1} \geq_{n+k+1} p_k$,
\item sets $\left\{F"\lft2(\lft1[(p_{k+1})_s\rgt1]\rgt2):  s \in
  \spli_{n+k+1}(p_{k+1})\right\}$ are pairwise disjoint and have
diameter $< 2^{-k}$
\end{enumerate}

Suppose that  $p_k$ is given. Note that $F"\lft1([(p_k)_s]\rgt1)$ is uncountable for every $s \in
p_k$. For $v \in \spli_{n+k+1}(p_k)$ choose pairwise different  reals $x_v \in
F"\lft1([(p_k)_v]\rgt1)$. It is not important now but will be relevant
in the sequel, that we can choose these reals ``effectively'' from a
fixed countable subset of $[p_k]$. Let $\ell >k $ be such that sequences $x_v \rest
\ell$ are also pairwise different. For every $v\in \spli_{n+k+1}(p_k)$ let $s_v \in
\spli(p_k)$ 
be such that for every $z \in [(p_k)_{s_v}]$, $F(z) \rest \ell=x_v
\rest \ell$. If $F$ is as in the second part of lemma \ref{general}
then we can find $s_v$ in $\spli_\ell(p_k)$. 
Define 
$$p_{k+1} = \bigcup \{(p_k)_{s_v}: v \in \spli_{n+k+1}(p_k)\}.$$
Observe that $q = \lim_k p_k$ has the required property. 
\end{proof}
 
Note that the above lemma shows that the reals added by $\EE$ are
minimal. Infinitely equal forcing from \cite{Mil81Som} or
\cite{Cor89Gen} does not have this property. 

\section{Iteration of $\EE$.}

Let $ \alpha \leq \omega_2$ be an ordinal and
suppose that $\EE_\alpha $ is a
countable support iteration of $\EE$ of length $\alpha $.
In other words, $p \in \EE_\alpha $ is
\begin{enumerate}
\item $p$ is a function and $\dom(p)=\alpha $,
\item $\supp(p)=\{\beta: p(\beta)\neq \emptyset\}$ is countable,
\item $\forall \beta<\alpha \ p\rest \beta \forces_{\EE_\beta }
  p(\beta) \in \EE$.
\end{enumerate}
For $F \in [\alpha]^{<\omega}$, $n \in \omega $, and $p,q \in \EE_\alpha $ define
$$q \geq_{F,n} p \iff q \geq p \ \&\ 
\forall \beta \in F \ q \rest \beta \forces_{\EE_\beta}
q(\beta)\geq_n p(\beta).$$

The following fact is well-known.
\begin{theorem}[\cite{GJS92}, \cite{Mil81Som}, \cite{BJbook}]\label{oldmil}
Suppose that $p \in \EE_\alpha $, $F \in [\alpha]^{<\omega}$, and $n
\in \omega $.
  \begin{enumerate}
  \item for every maximal antichain ${\mathcal A} \subseteq \EE_\alpha$,
    there exists $q \geq_{F,n} p$ such that $\{r \in
    {\mathcal A} : r \text{ is compatible with } q\}$ is finite.
  \item for every family of maximal antichains $\{{\mathcal A}_n: n
    \in \omega \}$
    there exists $q \geq p$ such that for every $n$, 
$\{r \in
    {\mathcal A}_n : r \text{ is compatible with } q\}$ is finite. 
\item $\V^{\EE_{\omega_2}} \thinks [\reals]^{<2^{\boldsymbol\aleph_0}}
  \subseteq \N$.
\item $\V^{\EE_{\omega_2}} \thinks \M \cap \V \text{ is cofinal in
    }\M$. $\qed$
\end{enumerate}
\end{theorem}

For $p \in \EE_{\alpha}$ let
$\cl(p)$ be the smallest set $w \subseteq \alpha$ 
such that $p$ can be evaluated using generic reals $\<\dot{g}_\beta :
\beta \in w\>$.
In other words, $\cl(p)$ consists of those $\beta< \alpha $ such that
the transitive closure of $p$ contains $\EE_\beta $ name for an
element of $\EE$. 
It is well-known \cite{shelahbook}
that
$\{p \in \EE_{\alpha}: \cl(p) \in [\alpha]^{\leq \omega}
    \}$ is dense in $\EE_{\alpha}$.

Suppose that $p \in \EE_\alpha $, $w=\cl(p)$ is countable and $\alpha_p=\ot\lft1(\cl(p)\rgt1)$. 
Let $\EE_w$ be the countable support iteration of $\EE$ with the
domain $w$. In other words, consider The countable support iteration
$\<{\mathcal P}_\beta, \dot{{\mathcal Q}}_\beta: \beta<\sup(w)\>$ such
that 
$$\forall \beta < \sup(w) \ \forces_{{\mathcal P}_\beta}
\dot{{\mathcal Q}}_\beta \simeq \left\{
  \begin{array}{ll}
\EE & \text{if } \beta\in w\\
\emptyset & \text{if } \beta \not\in w
  \end{array}\right. .$$
It is clear that $\EE_w \simeq \EE_{\alpha_p}$. 
Moreover, we can view condition $p$ as a member of $\EE_w$. 


For the rest of the section we will consider only the iteration of
$\EE$ of countable length $\alpha $ and show that $\EE_\alpha $ has
the same properties that $\EE$. 

Let $\alpha $ be a countable ordinal and $p \in \EE_\alpha $.
Define $\overline{p} \subseteq
\X^{\alpha} $ as follows:

$\<x_\beta:\beta<\alpha \> \in \overline{p}$ if
for every $\beta < \alpha $,
$$x_\beta \in \lft2[p(\beta)\lft1[\<x_\gamma:\gamma<\beta\>\rgt1]\rgt2].$$
Note that $p(\beta)[\<x_\gamma:\gamma<\beta\>]$ is the interpretation
of $p(\beta)$ using reals $\<x_\gamma:\gamma<\beta\>$ so may be
undefined if these reals are not sufficiently generic.

For a set $G \subseteq \X^\alpha $, $u \subseteq \alpha $, and
$x \in \X^u$
let $$(G)_x=\{y \in \X^{\alpha \setminus u}  : \exists z\in G \ z \rest
u=x \ \&\ z \rest (\alpha \setminus u)=y\},$$
and for $\beta \in \alpha $ let
$$(G)_\beta = \{x(\beta): x \in G\}.$$ 

We say that $p \in \EE_\alpha $ is good if 
\begin{enumerate}
\item $\overline{p}$ is compact,
\item for every $\beta<\alpha $
and $x \in \overline{p\rest \beta}$, $\overline{p[x]}=(p)_x$ and 
$\overline{p(\beta)[x]}=((p)_x)_\beta $.
\item $\overline{p}$ is homeomorphic to $\X^\alpha $ via a
  homeomorphism $h$ such that  for every $\beta < \alpha $ and $x \in
  \overline{p \rest \beta}$, $h \rest ((p)_x)_\beta$ is a
  homeomorphism between $((p)_x)_\beta$ and $\X$.
\end{enumerate}


\begin{lemma}
  $\{p \in \EE_\alpha : \overline{p} \text{ is good}\}$ is dense in
  $\EE_\alpha $.
\end{lemma}
\begin{proof}
{\sc Case 1}. $\alpha=\beta+1$. 

Fix $p \in \EE_\alpha $ and for $n \in \omega $ let ${\mathcal A}_n$
be a maximal antichain below $p\rest \beta $ such that 
\begin{enumerate}
\item $\forall r \in {\mathcal A}_n\ \overline{r}$ is compact.
\item $\forall r \in {\mathcal A}_n\ \exists t \subseteq \prod_{j<n}f(j) \ r \forces_{\EE_\beta}
  p(\beta) \rest n =t$.
\end{enumerate}
Fix a sequence $\<F_n: n \in \omega \>$ such that for $n \in \omega $, 
\begin{enumerate}
\item $F_n \in [\beta]^{<\omega}$,
\item $F_n \subseteq F_{n+1}$,
\item $\bigcup_n F_n = \beta $.
\end{enumerate}
By induction build a sequence $\<q_n: n \in \omega\>$ such that for $n
\in \omega $,
\begin{enumerate}
\item $\overline{q_n}$ is compact,
\item $q_{n+1} \geq_{F_n,n} q_n$,
\item $\exists {\mathcal A}_n'\in [{\mathcal A}_n]^{<\omega}\
  \overline{q_n} \subseteq \bigcup_{r \in {\mathcal A}_n'} \overline{r}$. 
\end{enumerate}
Let $q_\omega  = \lim_n q_n$. As in the proof of \ref{general} we show that
there exists a continuous function $F: \overline{q_\omega } \longrightarrow
\EE$ (encode elements of $\EE$ as reals) such that 
$$q_\omega  \forces_{\EE_\beta}
p(\beta)=F(\<\dot{g}_\gamma:\gamma<\beta\>).$$
Consider $q = q_\omega{}^\frown p(\beta)\geq p$.
Clearly,
$$\overline{q}=\{\<x,y\>: x \in \overline{q_\omega}, \ y \in
[F(x)]\}$$
is compact in $\X^\alpha $. Remaining requirements are met as well.

\bigskip

{\sc Case 2}. $ \alpha $ is limit.

Given $p \in \EE_\alpha $ fix sequences $\<F_n: n \in \omega\>$ and
$\<\alpha_n: n \in \omega\>$ such that 
\begin{enumerate}
\item $F_n \in [\alpha_n]^{<\omega}$,
\item $F_n \subseteq F_{n+1}$,
\item $\bigcup_n F_n = \alpha $,
\item $\sup_n \alpha_n=\alpha $.
\end{enumerate}
By induction build a sequence $\<q_n: n \in \omega\>$ such that for $n
\in \omega $,
\begin{enumerate}
\item $q_n \in \EE_\alpha $,
\item $\supp(q_n) \subseteq \alpha_n$,
\item $q_{n+1} \geq_{F_n,n} q_n$,
\item $q_n \rest \alpha_n \geq p \rest \alpha_n$,
\item $\overline{q_n \rest \alpha_n} $ is compact in $\X^{\alpha_n}$.
\end{enumerate}
Let $q = \lim_n q_n$.
Note that $\overline{q}=\bigcap_n \overline{q_n \rest \alpha_n} \times
\X^{\alpha \setminus \alpha_n}$ is as required.
\end{proof}

From now on we will always work with conditions $p$ such that
$\overline{p}$ is good.
We noticed earlier that for every condition $p \in \EE$, $[p]$ is canonically
isomorphic to $2^\omega $, in exactly the same way we can verify that
if $p \in \EE_\alpha $ and $\overline{p}$ is
good then $\overline{p}$ is isomorphic to $(2^\omega)^\alpha $.

As in the lemma \ref{general} we show that:
\begin{lemma}\label{old1}
  Suppose that $p \in \EE_\alpha $ and $p \forces_{\EE_\alpha} \dot{x}
  \in 2^\omega $.
Then there exists $q \geq p$ and a continuous function $F:
\overline{p} \longrightarrow 2^\omega $ such that 
$$q \forces_{\EE_\alpha } \dot{x} =F(\dot{\mathbf g}),$$
where $\dot{{\mathbf g}}=\<\dot{g}_\beta : \beta<\alpha \>$ is the sequence of generic reals.
\end{lemma}

\begin{lemma}\label{cruc}
  Let $p \in \EE_\alpha $ and suppose that $H \subseteq \overline{p}$ is a
  meager set in $\overline{p}$. For every $F \in [\alpha]^{<\omega}$
  and $n \in \omega $ there
  exists $q \geq_{F,n} p$ such that $\overline{q} \cap H=\emptyset$. 
\end{lemma}
\begin{proof}
As before, without loss of generality we can assume that $\alpha $ is
countable.

Induction on $\alpha $. 

{\sc Case 1}. $ \alpha =\beta+1$.

Suppose that $p \in \EE_\alpha $ and $H \subseteq \overline{p}
\subseteq \X^\beta \times \X$ is meager, and let $F \in
[\alpha]^{<\omega} $ and $n \in \omega $ be given.

Let 
$$H'=\left\{x \in \overline{p \rest \beta }: (H)_x \text{ is not meager in
  } \lft1[p(\beta)[x]\rgt1]=((\overline{p})_x)_\beta\right\}.$$
Using the fact that $\overline{p}$ is homeomorphic to
$(2^\omega)^\alpha $ via homeomorphism respecting vertical sections,
and 
by Kuratowski-Ulam theorem, we conclude that  $H'$ is a meager set in $\overline{p\rest
  \beta}$.

Recall the following classical lemma:
\begin{lemma}[\cite{MeasuCatInv}]\label{fremlin}
Suppose that $H \subseteq 2^\omega \times 2^\omega$ is a Borel set.
\begin{enumerate}
\item Assume $(H)_x$ is meager for all $x$. Then there
exists a sequence of Borel sets $\{G_n :n \in \omega\} \subseteq
2^\omega \times 2^\omega $ such that
\begin{enumerate}
  \item $(G_n)_x$ is a closed nowhere dense set for all $x \in
2^\omega$,
\item $H \subseteq \bigcup_{n \in \omega} G_n$.
\end{enumerate}
\end{enumerate}
  \end{lemma}

By the inductive hypothesis we can find $q^\star \geq_{F\cap \beta, n} p \rest \beta $ such
that $\overline{q^\star} \cap H' =\emptyset$.
By lemma \ref{spec} for every $x \in \overline{q^\star}$ there exists $q_x \geq_n
p(\beta)[x]$ such that $[q_x] \cap (H)_x = \emptyset$.
Moreover, by \ref{fremlin}, the mapping $x \mapsto q_x$ is can be
chosen to be Borel, and subsequently, by shrinking $q^\star$,
continuous. 
Let $q \in \EE_\alpha $ be defined such that 
$q \rest \beta =q^\star$ and $q^\star \forces_{\EE_\beta} q(\beta)=q_{\dot{g}_\beta}$. 
It is clear that $q$ has the required properties.

{\sc Case 2}. $\alpha $ is limit.

Fix sequences $\<F_n: n \in \omega\>$ and
$\<\alpha_n: n \in \omega\>$ such that 
\begin{enumerate}
\item $F_n \in [\alpha_n]^{<\omega}$,
\item $F_n \subseteq F_{n+1}$,
\item $\bigcup_n F_n = \alpha $,
\item $\sup_n \alpha_n=\alpha $.
\end{enumerate}
By induction build a sequence $\<q_n: n \in \omega\>$ such that for $n
\in \omega $,
\begin{enumerate}
\item $q_n \in \EE_\alpha $,
\item $\supp(q_n) \subseteq \alpha_n$,
\item $q_{n+1} \geq_{F_n,n} q_n$,
\item $q_n \rest \alpha_n \geq p \rest \alpha_n$,
\item $\overline{q_n \rest \alpha_n} \cap H_n =\emptyset$, where
$H_n=\left\{x \in \overline{q_n \rest \alpha_n}: (H)_x \text{ is not meager
  in } \overline{p[x]}\right\}$.
\end{enumerate}
As before (5) is possible by  Kuratowski-Ulam theorem.
Let $q = \lim_n q_n$. It is clear that $\overline{q} \cap  H =
\emptyset$. 
\end{proof}

The following lemma is an analog of lemma \ref{uff}.
\begin{lemma}\label{uff1}
  Suppose that $p \in \EE_\alpha$, $n \in \omega $ and $p
  \forces_{\EE_\alpha} \dot{x} \in 2^\omega
  $. Let $F : \overline{p} \longrightarrow 2^\omega $ be a continuous function
  such that $p \forces_{\EE_\alpha} \dot{x}=F(\dot{\mathbf g})$, 
where $\dot{{\mathbf g}}=\<\dot{g}_\beta : \beta<\alpha \>$ is the sequence of generic reals.
There exists $q \geq p$ 
such that exactly one of the following conditions hold:
\begin{enumerate}
\item $F \rest \overline{q} $ is constant,
\item there exists $\beta< \alpha $ such that $F \rest \overline{q
    \rest \beta}$ is one-to-one and for every $x \in \overline{q \rest
    \beta}$, $F \rest (\overline{q\rest \beta})_x$ is constant,
\item $F \rest \overline{q}$ is one-to-one.
\end{enumerate}
\end{lemma}
\begin{proof}
We have three  cases:

{\sc Case 1}. There exists $q \geq p$ such that $q \forces_{\EE_\alpha}
\dot{x} \in \V$. 
Without loss of generality we can assume that for some $x \in \V \cap
2^\omega $
$q
\forces_{\EE_\alpha} \dot{x}=x$. It follows that $F \rest
\overline{q}$ is constant.

\bigskip 

{\sc Case 2}. There exists $q \geq p$ such that $q \forces_{\EE_\alpha}
\exists \beta< \alpha \ \dot{x} \in \V^{\EE_\beta}$.
By shrinking $q$  we can assume that there exists a continuous
function $G: \overline{q \rest \beta} \longrightarrow 2^\omega $ such
that $q \forces_{\EE_\alpha} \dot{x}=G(\dot{\mathbf g}\rest
\beta)$. In particular,
for $x \in [q]$, $F(x)=G(x \rest \beta)$.
If $\beta $ was minimal then, using the argument below, we can also
assume that $G$ is one-to-one.

\bigskip

Suppose that $q \in \EE_\alpha $, $F \in [\alpha]^{<\omega}$, 
 and $n \in \omega $.
Without loss of generality we can assume that for every $\beta \in F$, $q
\rest \beta$ determines the value of $\spli_n(q(\beta))$ (up to
finitely many values).
Suppose that $\sigma:F \longrightarrow \omega^{<\omega}$ is a function
such that $\sigma(\beta) \in \spli_n(q(\beta))$ for $\beta \in F$.
Let $(q)_\sigma $ be the condition defined as
$$\forall \beta < \alpha \ (q)_\sigma\rest \beta  \forces_{\EE_\beta}
(q)_\sigma(\beta)=\left\{
  \begin{array}{ll}
q(\beta)& \text{if } \beta \not \in F\\
(q(\beta))_{\sigma(\beta)}& \text{if } \beta \in F
  \end{array}\right. . $$
Let $\Sigma_{F,n}$ be the finite set of all mappings $\sigma $
satisfying the requirements. 

\begin{lemma}
  Suppose that $F \in [\alpha]^{<\omega}$, $n \in \omega $ and 
$$p \forces_{\EE_\alpha} \dot{x} =F(\dot{{\mathbf g}}) \
  \& \ \forall \beta<\alpha \ \dot{x}\not\in \V^{\EE_\beta}.$$
There exists $q \geq_{F,n} p$ such that the sets
$\left\{F"(\overline{(q)_\sigma }): \sigma \in \Sigma_{F,n}\right\}$ are pairwise
disjoint.
\end{lemma}
\begin{proof}
 Induction on $|F|$ and $\alpha$.
If $F=\{\beta\}$ this is essentially lemma \ref{uff}.

Let $\{v_j: j<k^\star\}$ be an enumeration of $\spli_n(p(\beta))$.
For $v \in \spli_{n}(p)$ choose pairwise different  reals $x_v \in
F"\lft1(\overline{(p)_v}\rgt1)$. Note that this choice can be made
canonically from, for example, the countable dense let of leftmost
branches of subtrees of $p$.  Let $\ell >k $ be such that sequences $x_v \rest
\ell$ are also pairwise different. 
Define conditions $\<r_j: j \leq k^\star\>$, $\<q_j: j \leq  k^\star\>$ such
that for every $j \leq k^\star$,
\begin{enumerate}
\item $r_j \in \EE_\beta $,
\item $r_{j+1} \geq r_j$,
\item $r_j \forces_{\EE_\beta} q_j \geq (p)_{v_j} \rest
  [\beta,\alpha)$,
\item $\forall z \in \overline{r_j{}^\frown q_j}, \ F(z) \rest \ell =
  F(x_{v_j}) \rest \ell$.
\end{enumerate}
Let $q \rest \beta = q_{k^\star}$ and $q \rest [\beta, \alpha ) =
\bigcup_{j<k^\star} q_j$.

\bigskip

Suppose that $|F|=k+1$ and let $\beta=\max(F)$.

By the part already proved,
for each ${\mathbf x}=\<x_\gamma: \gamma<\beta\> \in \overline{p\rest
  \beta}$ find a condition $q_{\mathbf x} \geq_n p\rest [\beta, \alpha
)[{\mathbf x}]$ such that the sets $\left\{F"(\overline{(q_{\mathbf x})_s}):
s \in \spli_n(q_{\mathbf x})\right\}$ are pairwise disjoint.
Note that we can do it in such a way that the mapping ${\mathbf x}
\mapsto q_{\mathbf x}$ is continuous (As before we first choose
$q_{\mathbf x}$ in a Borel way, and then shrink $p \rest
\beta$ to make this mapping continuous). That defines a $\EE_\beta
$-name for an element of $\EE_{\beta, \alpha }$, which we call
$q^\star$.

Next, let $F'=F \setminus \{\beta\}$ and apply the inductive
hypothesis to find $q' \geq_{F', n} p \rest \beta $ such that 
$\left\{F"(\overline{(q')_\sigma }): \sigma \in \Sigma_{F',n}\right\}$ are
pairwise disjoint.
Let $q \in \EE_\alpha $ be defined as $q \rest \beta=q'$ and $q \rest
\beta \forces_{\EE_\beta} q \rest [\beta, \alpha )=q^\star$.

It is clear that $q$ is as required. 
\end{proof}

\bigskip

{\sc Case 3}. $p \forces_{\EE_\alpha} \forall \beta<\alpha \ \dot{x}
\not\in \V^{\EE_\beta}$.

Let $\<F_n: k \in \omega \>$ be an increasing sequence of finite sets
such that $\bigcup_n F_n =\alpha $.

By induction build a sequence of conditions $\<p_n: n \in \omega\>$
such that $p_0=p$ and for every $n$,
\begin{enumerate}
\item $p_{n+1} \geq_{F_n,n} p_n$,
\item sets $\left\{F"(\overline{(p_n)_\sigma }): \sigma \in
  \Sigma_{F_n,n}\right\}$ are pairwise disjoint.
\end{enumerate}
Let $q =\lim_n p_n$. 

Suppose that ${\mathbf x}=\<x_\beta: \beta<\alpha \> $ and ${\mathbf
  x}'=\<x'_\beta: \beta<\alpha \> $ are two distinct points in
$\overline{q}$.
Let $\beta $ be the first ordinal such that $x_\beta \neq x_\beta'$.
Let $n$ be so large that
$\beta \in F_n$ and there are two distinct $\sigma, \sigma' \in
\Sigma_{F_n,n}$ such that ${\mathbf x} \in \overline{(p_n)_\sigma }$
and ${\mathbf x}' \in \overline{(p_n)_{\sigma'} }$.
Since $F"(\overline{(p_n)_\sigma}) \cap
F"(\overline{(p_n)_{\sigma'}})=\emptyset$, 
it follows that $F({\mathbf x})\neq
F({\mathbf x}')$. 
\end{proof}

\section{A model where $\PM \subseteq \UN$.}

Let $\EE_{\omega_2}$ be the countable support iteration of $\EE$ of
length ${\boldsymbol\aleph_2}$. 
We will show that in $\V^{\EE_{\omega_2}}$, $\PM \subseteq \UN$.

By theorem \ref{oldmil}(2), $\V^{\EE_{\omega_2}} \thinks [\reals]^{<2^{\boldsymbol\aleph_0}}
  \subseteq \UN$,
thus we have to show that 
$$\V^{\EE_{\omega_2}} \thinks \PM \subseteq 
[\reals]^{<2^{\boldsymbol\aleph_0}}.$$

Suppose that $X \in \V^{\EE_{\omega_2}}$ is a set of reals of size
${\boldsymbol\aleph_2} $.
Let $\{\dot{x}_\alpha: \alpha<\omega_2\}$ be the set of names for
elements of $X$ such that $\forces_{\EE_{\omega_2}} \forall \alpha\neq
\beta \ \dot{x}_\alpha \neq \dot{x}_\beta $.
Apply lemma \ref{old1} and find for each $\alpha < \omega_2$ a set $w_\alpha \in [\omega_2]^{\leq
  \omega }$, a condition
$p_\alpha \in \EE_{w_\alpha}$,  and a continuous function $F_\alpha : \overline{p_\alpha}
\longrightarrow 2^\omega $ such that 
$$p_\alpha \forces_{\EE_{\omega_2}} \dot{x}_\alpha =
F_\alpha(\<\dot{g}_\beta: \beta\in w_\alpha\>).$$
We can assume that $w_\alpha $ is minimal. In other words, 
$$p_\alpha
\forces_{\EE_{\omega_2}} \forall \beta<\sup(w_\alpha) \ \dot{x}_\alpha
\not\in \V^{\EE_\beta}.$$
 In particular, without loss of generality we can
assume $F_\alpha $ is one-to-one, so it is a homeomorphism.
 
By thinning out we can assume that $\ot(w_\alpha)=\gamma$,
$F_\alpha=F$ and $\overline{p_\alpha}=\overline{p}$.
Moreover, since $\V \thinks \CH$, we can assume that $w_\alpha \cap
w_\beta = w^\star$ for $ \alpha \neq \beta $. Finally, without loss of
generality we can assume that $w^\star=\emptyset$.

Let $P=F"(\overline{p})$. Since $F$ is a homeomorphism, $P$ is perfect.
We will show that 
$X \cap P$ is not meager in $\V^{\EE_{\omega_2}}$ (relative to $P$).

Assume otherwise and let $H \subseteq P$ be a meager set such that for
some $p^\star \in \EE_{\omega_2}$, 
$p^\star \forces_{\EE_{\omega_2}} X \cap P
\subseteq H$. By
\ref{oldmil}(4) we can assume that $H \in \V$. 
Set $G=
(F)^{-1}(H)$ and notice that  $G$ is a meager subset of
$\overline{p}$.

Find $\alpha < \omega_2 $ such that $w_\alpha \cap
\cl(p^\star) = \emptyset$. 
By lemma \ref{cruc} there exists $q \geq p$, $q \in
\EE_{w_\alpha} \simeq \EE_\gamma$ such that 
$\overline{q} \cap G=\emptyset$.

Since $p^\star$
and $q$ are compatible let $r \geq p^\star, q$. 
It follows that
$$r \forces_{\EE_{\omega_2}} \dot{x}_\alpha =F_\alpha(\<\dot{q}_\beta: \beta\in w_\alpha
\>) \not \in H,$$
which finishes the proof.

{\bf Acknowledgments}:
The work was done while the first author  was spending  his sabbatical year at 
 the 
Rutgers University and the College of Staten
Island, CUNY, and  their support is gratefully acknowledged.


\end{document}